\documentclass{article}

\usepackage{arxiv}

\usepackage[utf8]{inputenc} 
\usepackage[T1]{fontenc}    
\usepackage{hyperref}       
\usepackage{url}            
\usepackage{booktabs}       
\usepackage{amsfonts}       
\usepackage{nicefrac}       
\usepackage{microtype}      
\usepackage{lipsum}		
\usepackage{graphicx}
\usepackage{natbib}
\usepackage{doi}

\usepackage{enumitem}		
\usepackage{indentfirst}		
\usepackage{color}			
\usepackage{graphicx}			
\usepackage{microtype} 
\usepackage{scalerel}[2016/12/29]
\usepackage{mathtools,amsfonts,amssymb,amsxtra,empheq,amsthm}
\usepackage{comment}
\usepackage{stackengine}

\newcommand{\disp}{\displaystyle}

\newcommand{\defi}{\coloneqq}

\newcommand{\bs}[1]{\boldsymbol{#1}}

\makeatletter
\newcommand*{\bigcdot}{}
\DeclareRobustCommand*{\bigcdot}{%
  \mathbin{\mathpalette\bigcdot@{}}%
}
\newcommand*{\bigcdot@scalefactor}{.5}
\newcommand*{\bigcdot@widthfactor}{1.5}
\newcommand*{\bigcdot@}[2]{%
  \sbox0{$#1\vcenter{}$}
  \sbox2{$#1\cdot\m@th$}%
  \hbox to \bigcdot@widthfactor\wd2{%
    \hfil
    \raise\ht0\hbox{%
      \scalebox{\bigcdot@scalefactor}{%
        \lower\ht0\hbox{$#1\bullet\m@th$}%
      }%
    }%
    \hfil
  }%
}
\makeatother

\newcommand{\N}{\mathbb{N}}

\newcommand{\R}{\mathbb{R}}

\newcommand{\varast}{\divideontimes}

\newtheorem{teorema}{Theorem}[section]
\newtheorem{proposicao}[teorema]{Proposition}
\newtheorem{definicao}[teorema]{Definition}

\newtheorem{exemplo}{Example}

{\theoremstyle{definition}

\newtheorem*{prova}{Prova}}

\title{A Carathéodory's Extension Theorem for Families Simpler Than an Algebras of Sets}

\author{ \href{https://orcid.org/0000-0002-6800-3289}{\includegraphics[scale=0.06]{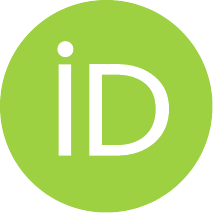}\hspace{1mm}Patrick Oliveira} \\
	Departamento de Matemática\\
	Universidade Federal de Minas Gerais\\
	\texttt{patrickoliveira@mat.dout.ufmg.br} \\
}

\renewcommand{\shorttitle}{\textit{arXiv} Template}

\hypersetup{
pdftitle={A Carathéodory Extension Theorem for Families Simpler Than an Algebras},
pdfsubject={q-mat.NC},
pdfauthor={Patrick Oliveira},
pdfkeywords={Measure Theory, Generation of Measures, Measure Extesion},
}

\begin{document}
\maketitle

\begin{abstract}

The Carathéodory's Extension Theorem is a powerful tool that allows us to generate a measure, over a sigma-algebra, from a pre-measure defined over an algebra of sets. However, although this result reduces our work to define a measure by only needing to define a pre-measure, it is not always easy to define the latter. The problem occurs when taking the smallest algebra that contains a family of targeted sets, it can be very complicated to consistently define the value of the pre-measure over its finite union of these sets - a union that is an element of the algebra. Thus, our objective in this article is to reproduce an extension theorem, just like the Carathéodory's Extension Theorem, but in the context of probability measures and replacing the need for a probability pre-measure defined over an algebra for now a quasi-measure defined over a refinement. The gain, then, is that the \textit{manual elaboration} of a quasi-measure is simpler than the elaboration of a pre-measure, since a refinement is a simpler structure than an algebra.
\end{abstract}

\keywords{Measure Theory \and Generation of Measures \and Measure Extesion}

\section{Introduction}

Imagine that we wish to construct a probability measure $\mathfrak{p}$ for a dynamical system whose possible states are given by a set $\Omega$ and whose dynamics are dictated by a family of simple events $\{\mathfrak{S} _\lambda\}_{\lambda} \subset \mathcal{P}(\Omega)$, where $\mathcal{P}(\Omega)$ is the set of partitions of $\Omega$. For example, we may be interested in the dynamic system given by the evolution of a communicable disease in a population. In this case, $\Omega$ can be given as all possible configurations in which this population is divided between healthy and sick individuals in contact through time. We can, for example, consider the family of simple events that dictate the evolution of this dynamic system as the sets $\mathfrak{S}_t$ given by the encounters between healthy and sick individuals in time $t \in \R$ . See that it is natural to define the probability of events $\mathfrak{S_t}$ occurring (just as it is natural for the complementary event $\mathcal{C}(\mathfrak{S}_t) \defi \Omega \setminus \mathfrak{S }_t$). Now, to actually define a probability over $\Omega$ capable of evaluating each event $\mathfrak{S}_t$, we need, using Carathéodory's Extension Theorem (\cite{bartle}), to at least define an algebra that contains $\{\mathfrak{ S}_t\}_{t \in \R}$ and a pre-measure $\mathfrak{p}^\ast$ defined over this algebra.


The problem that arises, then, is that although it is often easy to define the values  
$\mathfrak{p}^\ast(\mathfrak{S}_t)$ (or better, $\mathfrak{p}(\mathfrak{S}_t)$), it is not always intuitive how we could define, for example, the probabilities of the events $\cup_{i=1}^{n} \mathfrak{S}_{t_i}$ which, in turn, belong to any algebra that contains $\{\mathfrak{S}_t\}_{t \in \R}$. When the events $\{\mathfrak{S}_{t_i}\}_{i=1}^{n}$ are two-by-two disjoint, there is no problem, since we can define $\mathfrak{p}^\ast(\cup_{i=1}^{n} \mathfrak{S}_i) \defi 
\sum_{i=1}^{n} \mathfrak{p}^\ast(\mathfrak{S}_{t_i})$, however, we are not always guaranteed that the family of events $\{\mathfrak{S}_\lambda\}_{\lambda}$ is two-by-two disjoint in order to reproduce this intuitive definition. Even in the case of our example we can see that the events 
$\{\mathfrak{S}_t\}_{t \in \R}$ are not necessarily disjoint two by two, since, roughly speaking, a healthy individual can meeting another sick individual at different times.


\newpage

Therefore, it is interesting that there is a way to generate a probability without the headache of needing to define the value of a pre-measure on all the elements of an algebra. To this end, we will work with \textit{refinements} instead of algebras, which are a type of family of sets of $\Omega$ \textit{simpler} and, on these, we define a \textit{quasi-measure} that, just like a pre-measure in Carathéodory's Extension Theorem (\cite{bartle}), we can extend it to become a probability over $\Omega$.

\section{Basic Definitions}
\label{sec:basic_definitions}


To clarify the context in which we are working, let's be $\Omega$ and a collection $\mathfrak{S} \defi \{\mathfrak{S}_\lambda\}_{\lambda} \subset \mathcal{P} (\Omega)$, we say that the algebra of $\Omega$ generated by $\mathfrak{S}$, denoted $\mathcal{A}(\mathfrak{S})$, is given by the smallest algebra that contains $\mathfrak{S}$. It is not difficult to see that $\mathcal{A}(\mathfrak{S})$ is given by the elements: $\emptyset, \Omega$ and $\cup_{i=1}^{n} \mathfrak{E}_i $, where $\mathfrak{E}_i$ is equal to $\mathfrak{S}_\lambda$ or $\mathcal{C}(\mathfrak{S}_\lambda)$ for some $\mathfrak{S} _\lambda \in \mathfrak{S}$. Furthermore, we say that an application $\mathfrak{p}^\ast: \mathcal{A} \longrightarrow [0,1]$ is a probability pre-measure on an algebra $\mathcal{A}$ of $\Omega $ if $\mathfrak{p}^\ast$ satisfies the properties: $\mathfrak{p}^\ast(\emptyset) = 0$; $\mathfrak{p}^\ast(\mathfrak{A}) \geq 0$ for all $\mathfrak{A} \in \mathcal{A}$; $\mathfrak{p}^\ast(\cup_{i=1}^{n} \mathfrak{A}_i) = \sum_{i=1}^{n} \mathfrak{p}^\ast(\mathfrak{A}_i)$ for every finite union $\{\mathfrak{A}_i\}_{i=1}^{n} \subset \mathcal{A}$ disjoint two by two; and $\mathfrak{p}^\ast(\Omega) = 1$. All these definitions were based on the references \citep{bartle,halmos}.


Now we turn to the basic definitions.

\begin{definicao}
Let $\Omega$ and $\mathfrak{S} \subset \mathcal{P}(\Omega)$. We say that $\mathfrak{S}$ is a \textbf{coat of} $\bs{\Omega}$ if $\emptyset, \Omega \in \mathfrak{S}$.
\label{def_revestimento}
\end{definicao}


\begin{definicao}

Let $\mathfrak{S}$ a coat of $\Omega$. We define the \textbf{refinement of} $\bs{\mathfrak{S}}$, denoted by $\mathfrak{S}'$, by

\begin{equation*}
\mathfrak{S}' \defi \Big\{ X \cap Y, \; X \cap \mathcal{C}(Y) \; \Big| \; X,Y \in \mathfrak{S} \; \Big\} \cdot
\end{equation*}

\label{def_refinamento}
\end{definicao}



Note that in refinements, unlike in algebras, we are not interested in the elements given by the finite union of elements of $\mathfrak{S}$, but only in their two-by-two intersections (as well as in the intersections of the elements of $\mathfrak{ S}$ with its complements). Clearly, if $\mathfrak{S}$ is closed by complements, that is, if $\mathfrak{S}_\lambda \in \mathfrak{S}$ implies that $\mathcal{C}(\mathfrak{S} _\lambda) \in \mathfrak{S}$, then $\mathfrak{S}'$ is given just by $\mathfrak{S}' = \{X \cap Y \, | \, X,Y \in \mathfrak{S}\}$.


\begin{definicao}

Let $\mathfrak{S}$ be a coat of $\Omega$ and $\mathfrak{S}'$ be its refinement. We say that a map $\mathfrak{p}': \mathfrak{S}' \longrightarrow [0,1]$ is a \textbf{quasi-probability measure of} $\bs{\mathfrak{S}}$ (or a quasi-measure of probability of $\mathfrak{S}$) when, for all for $X,Y \in \mathfrak{S}$, the hypotheses apply: 

\begin{enumerate}

\item[]\emph{\textbf{(i)}} $\mathfrak{p}'(\emptyset) = 0$ and $\mathfrak{p}'(\Omega) =1$.

\item[]\emph{\textbf{(ii)}} $\mathfrak{p}'(X) = \mathfrak{p}' \big(X \cap Y\big) + \mathfrak{p}' \big(X \cap \mathcal{C}(Y) \big)$.

\item[]\emph{\textbf{(iii)}} Exists $W \in \mathfrak{S}'$ such that $\big( X \cap Y \big) \subset W$ and $\mathfrak{p}'(X \cap Y) = \mathfrak{p}'(W)$.

\item[]\emph{\textbf{(iv)}} Exists $Z \in \mathfrak{S}'$ such that $ \big( X \cap \mathcal{C}(Y) \big) \subset Z$ and $\mathfrak{p}'(X \cap \mathcal{C}(Y)) = \mathfrak{p}'(Z)$.

\item[]\emph{\textbf{(v)}} For all finite colection $\{S_n\}_{n=1}^{m} \subset \mathfrak{S}$ such that $X \subset \bigcup_{n=1}^{m} S_n$, we have $\mathfrak{p'}(X) \leq \sum_{n=1}^{m} \mathfrak{p}'(S_n)$.

\end{enumerate}

\label{def_quasi_medida}
\end{definicao}


Note that if $\mathfrak{p}^\ast: \mathcal{A}(\mathfrak{S}) \longrightarrow [0,1]$ is a probability pre-measure, then automatically $\mathfrak{p} ^\ast$, restricted to $\mathfrak{S}'$, is a quasi-probability measure of $\mathfrak{S}$ (and, of course, we have $\mathfrak{S} \subset \mathfrak{S} ' \subset \mathcal{A}(\mathfrak{S})$). Furthermore, in the Definition \ref{def_quasi_medida} we can restrict the existence of $W$ and $Z$ in items (iii) and (iv) to $\mathfrak{S}$ instead of $\mathfrak{S}'$ and, in fact, in the following results we will consider this more restricted case during the demonstrations. Finally, it is worth noting that in hypothesis (v) we simply need to verify the inequality over all finite and two-by-two disjoint collections (whenever it's possible).


In the same way as we can define an exterior measure from a pre-measure (\cite{bartle}), we can also define a exterior \textit{quasi-measure} from a quasi-measure.

\begin{definicao}

Let $\mathfrak{S}$ be a coat of $\Omega$. We define the \textbf{exterior quasi-measure of} $\bs{\mathfrak{S}}$ by $\mathfrak{p}^\varast : \mathcal{P}(\Omega) \longrightarrow \R \cup \{ +\infty\}$ given, for each $A \in \mathcal{P}(\Omega)$, by

\begin{equation*}
\mathfrak{p}^\varast(A) \defi \inf \Bigg\{ \sum_{n=1}^{m} \mathfrak{p}'(\mathfrak{S_n}) \; \Bigg| \; m \in \N \mbox{ and } A \subset \bigcup_{n =1}^{m} \mathfrak{S_n}, \mbox{ such that } \{\mathfrak{S_n}\}_{n = 1}^{m} \subset \mathfrak{S} \Bigg\},
\end{equation*}


\noindent where $\inf$ is taken over all finite collections $\{\mathfrak{S_n}\}_{n = 1}^{m} \subset \mathfrak{S}$ such that $A \subset \bigcup_{n = 1}^{m} \mathfrak{S_n}$.

\label{def_quasi_medida_exterior}
\end{definicao}


What becomes evident through the Definitions \ref{def_refinamento} and \ref{def_quasi_medida} is that, when we simplify the domain of our ``new pre-measure'' to something simpler, that is, when we migrate from an algebra to a refinement, we immediately complicate the hypotheses that this ``new pre-measure'' must satisfy. On the other hand, sometimes it is much more comfortable to check two more hypotheses in the definition of a quasi-measure than to consistently define a pre-measure over an entire algebra.

\newpage

\section{Main Results}

\begin{proposicao}

Let $\mathfrak{S}$ be a coat of $\Omega$ and $\mathfrak{p}': \mathfrak{S}' \longrightarrow [0,1]$ be a quasi-measure. Then $\mathfrak{p}^\varast: \mathcal{P}(\Omega) \longrightarrow [0,1]$ satisfies the following properties:

\begin{enumerate}

\item[]\emph{\textbf{(i)}} $\mathfrak{p}^\varast(\emptyset) = 0$ and $\mathfrak{p}^\varast(\Omega) = 1$.

\item[]\emph{\textbf{(ii)}} $\mathfrak{p}^\varast(A) \geq 0$ for all $A \subset \Omega$.

\item[]\emph{\textbf{(iii)}} $\mathfrak{p}^\varast(A) \leq \mathfrak{p}^\varast(B)$ for all $A,B \subset \Omega$ such that $A \subset B$.

\item[]\emph{\textbf{(iv)}} $\mathfrak{p}^\varast(X) = \mathfrak{p}'(X)$ for all $X \in \mathfrak{S}$.

\item[]\emph{\textbf{(v)}} For all family of sets $\{A_n\}_{n = 1}^{m} \subset \mathcal{P}(\Omega)$, we have $$\mathfrak{p}^\varast \left( \bigcup_{n =1}^{m} A_n \right) \leq \disp\sum_{n=1}^{m} \mathfrak{p}^\varast(A_n).$$

\end{enumerate}

\label{prop_quasi_medida_externa}
\end{proposicao}

\begin{prova}

Item (ii) follows immediately from the Definition \ref{def_quasi_medida_exterior}, since $\mathfrak{p}'$ is non-negative. For item (i) let's see that, $\emptyset \in \mathfrak{S}$ and, therefore, we have $\mathfrak{p}^\varast(\emptyset) \leq \mathfrak{p}'(\emptyset) = 0$, where $\mathfrak{p}'(\emptyset) = 0$ since $\mathfrak{p}'$ is a quasi-measure. For item (iii) we just need to note that every coverage of $B \subset \Omega$ by collections of sets $\{S_n\}_{n = 1}^{m} \subset \mathfrak{S}$ is also coverage of $A \subset \Omega$, given that $A \subset B$.


As for item (iv), first let us see that, for each $X \in \mathfrak{S}$, we have $\mathfrak{p}^\varast(X) \leq \mathfrak{p}'(X)$. To show the opposite inequality, let $\{S_n\}_{n \in\N} \subset \mathfrak{S}$ be an arbitrary collection such that $X \subset \bigcup_{n \in \N} S_n$ . Therefore, as $\mathfrak{p}'$ is a quasi-measure, it follows that $\mathfrak{p}'(X) \leq \sum_{n=1}^{\infty} \mathfrak{p}'(S_n )$. Due to the arbitrary nature of the collection $\{S_n\}_{n \in \N}$, we have that $\mathfrak{p}'(X) \leq \mathfrak{p}^\varast(X)$.


For item (v), let $\{A_n\}_{n =1}^{m}$ be a collection of subsets of $\Omega$ and $\varepsilon > 0$. Note that for each $n \in \{1, \cdots, m\}$, there are indices $m_n$ and coverages $\{S_i^n\}_{i \in \{1, \cdots, m_n\}} \subset \mathfrak{S}$ such that $A_n \subset \bigcup_{i = 1}^{m_n} S_i^n$ and $\sum_{i=1}^{m_n} \mathfrak{p}' (S_i^n) \leq \mathfrak{p}^\varast(A_n) + \varepsilon/2^n$. Therefore, the finite collection $\{S_i^n\}_{i,n \in \{1, \cdots, m_n \} } \subset \mathfrak{ S}$ covers $\bigcup_{n =1}^{m} A_n$, where we have $$\mathfrak{p}^\varast\left( \bigcup_{n =1}^{m} A_n \right) \; \leq \; \sum_{n=1}^{m} \sum_{i=1}^{m_n} \mathfrak{p}'(S_i^n) \; \leq \; \sum_{n=1}^{m} \Big( \mathfrak{p}^\varast (A_n) + \varepsilon/2^n \Big) \; \leq \; \varepsilon + \sum_{n=1}^{m} \mathfrak{p}^\varast(A_n).$$


For the generality of $\varepsilon > 0$ the result follows. Finally, from item (iv) it follows that $\mathfrak{p}^\varast(\Omega) = \mathfrak{p}'(\Omega) = 1$. \hfill $\square$
\end{prova}


\vspace{0.25cm}


\begin{teorema}

Let $\mathfrak{S}$ be a coat of $\Omega$ and $\mathfrak{p}': \mathfrak{S}' \longrightarrow [0,1]$ be a quasi-measure of $\mathfrak{S}$. Then, the restriction of $\mathfrak{p}^\varast$ to the algebra $\mathcal{A}(\mathfrak{S})$, $\mathfrak{p}^\varast: \mathcal{A}(\mathfrak {S}) \longrightarrow [0,1]$, is a probability pre-measure.
\label{teo_pequeno_caratheodory}
\end{teorema}

\vspace{0.15cm}

\begin{prova}

Define $\mathfrak{A} \subset \mathcal{P}(\Omega)$ as follows, $$\mathfrak{A} \defi \Big\{W \subset \Omega \; \Big| \; \mbox{For all } A \subset \Omega \mbox{ we have } \mathfrak{p}^\varast(A) = \mathfrak{p}^\varast(A \cap W) + \mathfrak{p}^\varast(A \cap \mathcal{C}(W)) \Big\}.$$


First let us see that $\mathfrak{A}$ is, in fact, an algebra of $\Omega$. Clearly $\emptyset, \Omega \in \mathfrak{A}$. Thus, let us verify that $\mathfrak{A}$ is closed by intersection and complementation. In effect, let $W,Z \in \mathfrak{A}$. Therefore, for all $A \subset \Omega$ the following identities hold:

\begin{eqnarray*}
\mathfrak{p}^\varast(A) & \overset{(a)}{=} & \mathfrak{p}^\varast(A \cap W) + \mathfrak{p}^\varast(A \cap \mathcal{C}(W)), \\
\mathfrak{p}^\varast(A \cap W) & \overset{(b)}{=} & \mathfrak{p}^\varast \Big((A \cap W) \cap Z \Big) + \mathfrak{p}^\varast \Big( (A \cap W) \cap \mathcal{C}(Z)\Big), \\
\mathfrak{p}^\varast(A \cap \mathcal{C}(W \cap Z)) & \overset{(c)}{=} & \mathfrak{p}^\varast \Bigg( \Big(A \cap \mathcal{C}(W \cap Z) \Big) \cap W \Bigg) +  \mathfrak{p}^\varast \Bigg( \Big(A \cap \mathcal{C}(W \cap Z) \Big) \cap \mathcal{C}(W) \Bigg) \\
& = & \mathfrak{p}^\varast \Bigg( \Big(A \cap \big(\mathcal{C}(W) \cup \mathcal{C} (Z) \big) \Big) \cap W \Bigg) +  \mathfrak{p}^\varast \Bigg( \Big(A \cap \big( \mathcal{C}(W) \cup \mathcal{C}(Z) \big) \Big) \cap \mathcal{C}(W) \Bigg) \\
& = & \mathfrak{p}^\varast \Bigg( \Big( \big(A \cap \mathcal{C}(W)\big) \cup \big( A \cap \mathcal{C} (Z) \big) \Big) \cap W \Bigg) \\
& & \; + \; \mathfrak{p}^\varast \Bigg( \Big( \big(A \cap \mathcal{C}(W)\big) \cup \big( A \cap \mathcal{C} (Z) \big) \Big) \cap \mathcal{C}(W) \Bigg) \\
& = & \mathfrak{p}^\varast(A \cap \mathcal{C}(Z) \cap W) + \mathfrak{p}^\varast(A \cap \mathcal{C}(W)), 
\end{eqnarray*} 


\noindent where (a) is true because $W \in \mathfrak{A}$, (b) is true because $Z \in \mathfrak{A}$ and $A \cap W \subset \Omega$, and (c) is true because $W \in \mathfrak{A}$ and $A \cap (\mathcal{C}(W \cap Z)) \subset \Omega$. Rewriting the last identity, we gain that 
\begin{equation*}
\mathfrak{p}^\varast(A \cap W \cap \mathcal{C}(Z)) \; = \; \mathfrak{p}^\varast(A \cap \mathcal{C}(W \cap Z)) - \mathfrak{p}^\varast(A \cap \mathcal{C}(W)).
\end{equation*}


Combining the previous equations, we obtain the identity
\begin{eqnarray*}
\mathfrak{p}^\varast(A) & = & \mathfrak{p}^\varast(A \cap W) + \mathfrak{p}^\varast(A \cap \mathcal{C}(W)) \\
& = & \mathfrak{p}^\varast (A \cap W \cap Z ) + \mathfrak{p}^\varast ( A \cap W \cap \mathcal{C}(Z)) + \mathfrak{p}^\varast(A \cap \mathcal{C}(W)) \\
& = & \mathfrak{p}^\varast (A \cap W \cap Z) + \Big( \mathfrak{p}^\varast(A \cap \mathcal{C}(W \cap Z)) - \mathfrak{p}^\varast(A  \cap \mathcal{C}(W)) \Big) + \mathfrak{p}^\varast(A \cap \mathcal{C}(W)) \\
& = & \mathfrak{p}^\varast (A \cap W \cap Z) + \mathfrak{p}^\varast(A \cap \mathcal{C}(W \cap Z)),
\end{eqnarray*}


\noindent which shows that $W \cap Z \in \mathfrak{A}$. Therefore $\mathfrak{A}$ is closed for finite intersections. To verify that $\mathfrak{A}$ is closed by complementation, let $W \in \mathfrak{A}$ and $A \subset \Omega$ be. Therefore, $$\mathfrak{p}^\varast(A) = \mathfrak{p}^\varast(A \cap W) + \mathfrak{p}^\varast(A \cap \mathcal{C}(W )) = \mathfrak{p}^\varast(A \cap \mathcal{C}(W)) + \mathfrak{p}^\varast(A \cap \mathcal{C}(\mathcal{C}(W )) ),$$
\noindent where $\mathcal{C}(W) \in \mathfrak{A}$. We thus prove that $\mathfrak{A}$ is an algebra of $\Omega$.


Now, let's see that $\mathfrak{A}$ contains $\mathfrak{S}$. In this spirit, we just need to check that, for all $X \in \mathfrak{S}$ and all $A \subset \Omega$, we have $\mathfrak{p}^\varast(A) = \mathfrak{p}^\varast(A \cap X) + \mathfrak{p}^\varast(A \cap \mathcal{C}(X))$. By item (v) of Proposition \ref{prop_quasi_medida_externa}, as $A = (A \cap X) \cup (A \cap \mathcal{C}(X))$, we have $\mathfrak{p}^\varast (A) \leq \mathfrak{p}^\varast(A \cap X) + \mathfrak{p}^\varast(A \cap \mathcal{C}(X))$. In order to verify the opposite inequality, let $\{S_n\}_{n = 1}^{m} \subset \mathfrak{S}$ be an arbitrary collection such that $A \subset \bigcup_{n =1 }^ {m} S_n$. Note that, given $X \in \mathfrak{S}$, as $\mathfrak{p}'$ is quasi-measure, then for each $n \in \{1, \cdots, m\}$ the identity holds $\mathfrak{p}'(S_n) = \mathfrak{p}'(S_n \cap X) + \mathfrak{p}'(S_n \cap \mathcal{C}(X))$. In this way, we have
\begin{equation*}
\sum_{n=1}^{m} \mathfrak{p}'(S_n) \; = \; \sum_{n=1}^{m} \mathfrak{p}'(S_n \cap X) + \mathfrak{p}'(S_n \cap \mathcal{C}(X)) \; = \; \sum_{n=1}^{m} \mathfrak{p}'(Z_n) + \mathfrak{p}'(W_n) \; \overset{(\ast)}{\geq} \;  \mathfrak{p}^\varast(A \cap X) + \mathfrak{p}^\varast(A \cap \mathcal{C}(X)),
\end{equation*}

\noindent where, for each $n \in \{1, \cdots,m\}$, $Z_n, W_n \in \mathfrak{S}$ are such that $\big(S_n \cap X \big) \subset Z_n$, $\big(S_n \cap \mathcal{C}(X) \big) \subset W_n$ and, also, $\mathfrak{p}'(S_n \cap X) = \mathfrak{p}' (Z_n)$ and $\mathfrak{p}'(S_n \cap \mathcal{C}(X)) = \mathfrak{p}'(W_n)$. Since the collections $\{Z_n\}_{n =1}^{m}, \{W_n\}_{n =1}^{m} \subset \mathfrak{S}$ respectively cover $\big ( A \cap X \big)$ and $ \big( A \cap \mathcal{C}(X) \big)$, follows the inequality $(\ast)$.


Based on the last inequality, and the generality of the coverage $\{S_n\}_{n \in \N}$, we conclude that the inequality $\mathfrak{p}^\varast(A) \geq \mathfrak{p}^\varast(A \cap X) + \mathfrak{p}^\varast(A \cap \mathcal{C}(X))$ holds. Therefore $X \in \mathfrak{A}$, which implies that $\mathfrak{S} \subset \mathfrak{A}$.

\vspace{0.15cm}


Since $\mathfrak{A}$ is algebra and $\mathcal{A}(\mathfrak{S})$ is the smallest algebra that contains $\mathfrak{S}$, then $\mathcal{A}(\mathfrak{ S}) \subset \mathfrak{A}$. Finally, to $\mathfrak{p}^\varast : \mathcal{A}(\mathfrak{S}) \longrightarrow [0,1]$ let's see that it is the pre-measure probability of $\Omega$. In effect, $\mathfrak{p}^\varast(\emptyset) = 0$ and $\mathfrak{p}^\varast(\Omega) = 1$ for item (i) of Proposition \ref{prop_quasi_medida_externa} and, for every finite and disjoint collection $\{E_i\}_{i=1}^{n} \subset \mathcal{A}(\mathfrak{S})$, with $n \in \N$, we have
\begin{eqnarray*}
\mathfrak{p}^\varast \left( \bigcup_{i=1}^{n} E_n \right) & \overset{(a)}{=} & \mathfrak{p}^\varast \left( \Bigg(\bigcup_{i=1}^{n} E_n\Bigg) \cap E_1 \right) +  \mathfrak{p}^\varast \left( \Bigg(\bigcup_{i=1}^{n} E_n\Bigg) \cap \mathcal{C}(E_1) \right)\\
& \overset{(b)}{=} & \mathfrak{p}^\varast \left(E_1\right) +  \mathfrak{p}^\varast \left( \bigcup_{i=2}^{n} E_n \right)\\
& = & \mathfrak{p}^\varast \left(E_1\right) +  \mathfrak{p}^\varast \left( \Bigg(\bigcup_{i=2}^{n} E_n\Bigg) \cap E_2 \right) +  \mathfrak{p}^\varast \left( \Bigg(\bigcup_{i=2}^{n} E_n\Bigg) \cap \mathcal{C}(E_2) \right)\\
& = & \mathfrak{p}^\varast \left(E_1\right)  + \mathfrak{p}^\varast \left(E_2\right)  + \mathfrak{p}^\varast \left( \bigcup_{i=3}^{n} E_n \right)\\
& \vdots & \\
& = & \sum_{i=1}^{n} \mathfrak{p}^\varast \left(E_i\right)  + \mathfrak{p}^\varast \left( \left(\bigcup_{i=n}^{n} E_i \right) \cap \mathcal{C}(E_n) \right)\\
& = & \sum_{i=1}^{n} \mathfrak{p}^\varast \left(E_i\right),
\end{eqnarray*}


\noindent where (a) is true because $E_1 \in \mathcal{A}(\mathfrak{S}) \subset \mathfrak{A}$ and (b) is true because the collection $\{E_i\}_{i =1}^{n}$ is disjoint two by two. The remaining identities are justified by repeating the algorithm of the first two lines. \hfill $\square$
\end{prova}


\newpage



As a natural consequence of the Theorem \ref{teo_pequeno_caratheodory}, using the Caratheodóry's Extension Theorem, it follows that every quasi-measure of probability $\mathfrak{p}'$ defined over a refinement $\mathfrak{S}'$ of a coat $\mathfrak{S}$ can be extended to a probability measure $\mathfrak{p}$ over the sigma-algebra generated by $\mathfrak{S}$. 

In the following result, we give alternative conditions for a map defined over a refinement $\mathfrak{S}'$ to be a quasi-measure.

\begin{proposicao}

Let $\mathfrak{S}$ be a coat of $\Omega$ and an application $\mathfrak{p}'': \mathfrak{S}' \longrightarrow [0,1]$, where $\mathfrak{S} '$ is refinement of $\mathfrak{S}$. If $\mathfrak{p}''$ satisfies, for every pair $X,Y \in \mathfrak{S}$, the following properties

\begin{enumerate}

\item[]\emph{\textbf{(i)}} $\mathfrak{p}''(\emptyset) = 0$ and $\mathfrak{p}''(\Omega) = 1$.

\item[]\emph{\textbf{(ii)}} $\mathfrak{p}''(X) \leq \mathfrak{p}´''(Y)$, when $X \subset Y$.

\item[]\emph{\textbf{(iii)}} $\mathfrak{p}''(X ) = \mathfrak{p}''(X \cap Y) + \mathfrak{p}'' (X \cap \mathcal{C}(Y))$.

\item[]\emph{\textbf{(iv)}} Exists $K,W \in \mathfrak{S}$ such that $K \subset \big( X \cap Y \big) \subset W$ and $\mathfrak{p}''(K) = \mathfrak{p}''( X \cap Y ) = \mathfrak{p}''( W)$.

\item[]\emph{\textbf{(v)}} Exists $Z \in \mathfrak{S}$ such that $\big( X \cap \mathcal{C}(Y) \big) \subset Z$ and $\mathfrak{p}''( X \cap \mathcal{C}(Y) ) = \mathfrak{p}''(Z)$.

\end{enumerate} 

Then $\mathfrak{p}''$ is a quasi-measure of probability of $\mathfrak{S}$.

\label{prop_quasi_medida_de_probabilidade}
\end{proposicao}

\begin{prova}

To verify that $\mathfrak{p}''$ is a quasi-measure of $\mathfrak{S}'$, we just need to check that, for every collection $\{S_n\}_{n = 1}^{m} \subset \mathfrak{S}$, with $m \in \N$, and $X \in \mathfrak{S}$ such that $X \subset \bigcup_{n = 1}^{m} S_n$, then we must have $\mathfrak{p}''(X) \leq \sum_{n=1}^{m} \mathfrak{p}''(S_n)$. In this spirit, we will show, by induction, that this statement is true for every finite collection $\{S_n\}_{n=1}^{m} \subset \mathfrak{S}$, $m \in \N$. Thus, for each $m \in \N$, we define the statement $\mathcal{A}(m)$ by
\begin{equation*}
\mathcal{A}(m) \; \defi \; \Bigg\{ \mbox{For all }  \{S_n\}_{n=1}^{m} \subset \mathfrak{S} \mbox{ e } X \in \mathfrak{S} \mbox{ such that } X \subset \bigcup_{n=1}^{m} S_n, \mbox{ then } \mathfrak{p}''(X) \leq \sum_{n=1}^{m} \mathfrak{p}''(S_n) \Bigg\} \cdot
\end{equation*}


Let's see that $\mathcal{A}(1)$ is true. In fact, by property (i) of the application $\mathfrak{p}''$, for every unitary collection $\{S_n\}_{n=1}^{1} \subset \mathfrak{S}$ and every $X \in \mathfrak{S}$ such that $X \subset \bigcup_{n=1}^{1} S_n$, we have $\mathfrak{p}''(X) \leq \mathfrak{p}''(S_1) = \sum_{n=1}^{1} \mathfrak{p}''(S_n)$.


Now suppose that $\mathcal{A}(m)$ is true for some $m \in \N$ and see that $\mathcal{A}(m+1)$ is also true. In effect, let $\{S_n\}_{n =1}^{m+1} \subset \mathfrak{S}$ and $X \in \mathfrak{S}$ be such that $X \subset \bigcup_{ n=1}^{m+1} S_n$. Therefore, given $S_1,S_2 \in \mathfrak{S}$, there are $K,W \in \mathfrak{S}$ such that $K \subset \big( S_1 \cup S_2 \big) \subset W$ and $\mathfrak{p}''(K) = \mathfrak{p}''(S_1 \cup S_2) = \mathfrak{p}''(W)$. Therefore, it is valid that $X \subset W \cup \bigcup_{n=3}^{m+1} S_n$, where the finite collection $\{W\} \cup \{S_n\}_{n=3 }^{m+1} \subset \mathfrak{S}$ has $m$ elements and, therefore, as $\mathcal{A}(m)$ is true by hypothesis, we have that $$\mathfrak{p} ''(X) \leq \mathfrak{p}''(W)+ \sum_{n=3}^{m+1} \mathfrak{p}''(S_n) = \mathfrak{p}''( K) + \sum_{n=3}^{m+1} \mathfrak{p}''(S_n).$$


Furthermore, we have the following identity
$$\mathfrak{p}''(K) = \mathfrak{p}''(K \cap S_1) + \mathfrak{p}''(K \cap \mathcal{C}(S_1)).$$


As $K \subset \big(S_1 \cup S_2 \big)$, so $\big(K \cap \mathcal{C}(S_1) \big) \subset S_2$. Furthermore, since $(K \cap S_1) \subset S_1$, we have
\begin{eqnarray*}
\mathfrak{p}''(X) & \leq & \mathfrak{p}''(K) + \sum_{n=3}^{m+1} \mathfrak{p}''(S_n) \; = \; \mathfrak{p}''(K \cap S_1) + \mathfrak{p}''(K \cap \mathcal{C}(S_1)) + \sum_{n=3}^{m+1} \mathfrak{p}''(S_n)\\ 
& \leq & \mathfrak{p}''(S_1) + \mathfrak{p}''(S_2) + \sum_{n=3}^{m+1} \mathfrak{p}''(S_n),
\end{eqnarray*}


\noindent which shows that the statement $\mathcal{A}(m+1)$ is true. \hfill $\square$
\end{prova}


\vspace{0.15cm}


\begin{exemplo}

Let $\Omega = \R^+$ and the coat $\mathfrak{S} = \{\emptyset, \R^+, [a,b] \, | \, a \leq b, \, a,b \in \R^+\}$. Therefore, by \emph{Definition \ref{def_refinamento}}, we gain $\mathfrak{S}' = \{\emptyset, \R^+, [a,b], [a,b), (a,b] , [u,a) \cup (b,v] \, : \, u \leq a \leq b \leq v, \, a,b,u,v \in \R^+\}$. Define $ \mathfrak{p}':\mathfrak{S}' \longrightarrow [0,1]$ by $\mathfrak{p}'(\emptyset) = 0$, $\mathfrak{p}'(\R^+) = 1$, $\mathfrak{p}'([a,b]) = \mathfrak{p}'([a,b)) = \mathfrak{p}'((a,b]) = e^{ -a} - e^{-b}$ and $\mathfrak{p}'([u,a) \cup (b,v]) = e^{-u} - e^{-a} + e^ {-b} - e^{-v}$ for all $a,b,u,v \in \R^+$. Let's see that $\mathfrak{p}'$ is a quasi-probability measure. 

Indeed, the hypotheses (i), (iii) and (iv) from \emph{Definition \ref{def_quasi_medida}} are trivially satisfied. The hypothesis (v) can be easily verified since, taking a finite collection $\{S_n\}_{n=1}^{m}$ \textbf{disjoint two by two} such that $X \subset \cup_{n=1}^{m} S_n$, with $X \in \mathfrak{S}$, as $S_n$ and $X$ are connected for each $n$, then there is necessarily  $n_0 \in \{1, \cdots, m\}$ such that $X \subset S_{n_0}$, where the inequality follows naturally. Finally, to verify hypothesis (ii), it sufices to analyze the non-trivial case in which $X = [u,v]$, $Y = [a,b]$ and $X \cap \mathcal{C}( Y) = [u,a) \cup (b,v]$, where we obtain $\mathfrak{p}'(X) = \mathfrak{p}'(X \cap Y) + \mathfrak{p}'( X \cap \mathcal{C}(Y))$.
\end{exemplo}

\newpage

\bibliographystyle{unsrtnat}
\bibliography{references}  






\end{document}